\documentclass{article}
\PassOptionsToPackage{numbers, compress}{natbib}
\usepackage[preprint]{neurips_2024}
\usepackage[utf8]{inputenc} % allow utf-8 input
\usepackage[T1]{fontenc}    % use 8-bit T1 fonts
\usepackage{hyperref}       % hyperlinks
\usepackage{url}            % simple URL typesetting
\usepackage{booktabs}       % professional-quality tables
\usepackage{amsfonts}       % blackboard math symbols
\usepackage{nicefrac}       % compact symbols for 1/2, etc.
\usepackage{microtype}      % microtypography
\usepackage{xcolor}         % colors
\usepackage{amsmath}
\usepackage{amssymb}
\usepackage{amsthm}
\usepackage{multirow}
\usepackage{mathtools}
\usepackage{cleveref}
\crefname{figure}{figure}{figures}
\usepackage{caption}
\usepackage{graphicx}
\usepackage{subcaption}
\usepackage[]{algorithm2e}

\newtheorem{theorem}{Theorem}
\newtheorem{remark}{Remark}
\newtheorem{example}{Example}

\newcommand{\wave}[1]{\ensuremath{\widetilde{#1}}}

\title{Mixed Newton Method for \\ Optimization in Complex Spaces}

\author{%
  Nikita Yudin \\
  MIPT,\,\,\,HSE~University,\,\,\,Innopolis~University,\,\,\,IITP~RAS,\,\,\,ISP~RAS,\,\,\,FRC~CSC~RAS\\
  \texttt{iudin.ne@phystech.edu}\\
  \And
  Roland Hildebrand \\
  Disko Laboratory, MIPT \\
  \texttt{khildebrand.r@mipt.ru} \\
  AI research lab, Skolkovo Institute of Science and Technology \\
  \texttt{R.Hildebrand@skoltech.ru} \\
  \AND
  Sergey Bakhurin \\
  MIPT \\
  \texttt{bakhurin.sa@mipt.ru} \\
  \And
  Alexander Degtyarev \\
  MIPT \\
  \texttt{degtyarev.aa@phystech.edu} \\
  \And
  Anna Lisachenko \\
  MIPT \\
  \texttt{lisachenko.am@phystech.edu} \\
  \And
  Ilya Kuruzov \\
  Disko Laboratory, MIPT\\
  IITP RAS,\,\,\,IAI MSU \\
  \texttt{kuruzov.ia@phystech.edu} \\
  \And
  Andrei Semenov \\
  MIPT,\,\,\,Yandex \\
  \texttt{semenov.a@phystech.edu}
  \\
  \And
  Mohammad Alkousa \\
  MIPT \\
  \texttt{mohammad.alkousa@phystech.edu}
}

\begin{document}

\maketitle

\begin{abstract}
  In this paper, we modify and apply the recently introduced Mixed Newton Method, which is originally designed for minimizing real-valued functions of complex variables, to the minimization of real-valued functions of real variables by extending the functions to complex space. We show that arbitrary regularizations preserve the favorable local convergence properties of the method, and construct a special type of regularization used to prevent convergence to complex minima. We compare several variants of the method applied to training neural networks with real and complex parameters.
\end{abstract}

\section{Introduction}

In this paper, we study the properties of the regularized mixed Newton method (RMNM). It is a development of the mixed Newton method (MNM) which was investigated in~\cite{Bakhurin2024}. We briefly recall the formulation and properties of the MNM and motivate its extension by regularization.

The MNM is used for the minimization of functions of the form
\[ f(z) = \sum_{j = 0}^{m - 1} |g_j(z)|^2,
\]
where the functions $g_j$ are holomorphic in the complex variable $z \in \mathbb C^n$. It is an iterative method starting at a point $z^0 \in \mathbb C^n$ and generating iterates according to the formula
\begin{equation} \label{MNMiterate}
    z^{k+1} = z^k - \left( \frac{\partial^2f(z^k)}{\partial \bar z\partial z} \right)^{-1} \frac{\partial f(z^k)}{\partial\bar z},
\end{equation}
where the derivatives are defined as Wirtinger derivatives and evaluated at the current point $z^k$.

The main advantages of the method are summarized in the following Theorem~\cite{Bakhurin2024}.

\begin{theorem} \label{thm:main}
The mixed derivative used in the MNM can be computed by the formula
\[ \frac{\partial^2f}{\partial \bar z\partial z} = \sum_{j=0}^{m - 1}\overline{\frac{dg_j}{dz}}\frac{dg_j}{dz}^\top
\]
and is hence positive semi-definite.

Let $\hat z$ be a critical point of the function $f$, i.e., $\frac{\partial f(\hat z)}{\partial\bar z} = 0_n$. In the non-degenerate case (when the full Hessian $\frac{\partial^2f(\hat z)}{\partial(z,\bar z)^2}$ is invertible) the iterates behave for $z$ near $\hat z$ asymptotically as if driven by a linear dynamical system, $z^{k+1} = \hat z + L(z^k-\hat z) + O(\|z^k-\hat z\|^2)$, where $L$ is an $\mathbb R$-linear operator. If $\hat z$ is a local minimum of $f$, then $L$ is contracting. If $\hat z$ is a saddle point of $f$, then $L$ has repulsive directions.
\end{theorem}

This means that the minima of $f$ are surrounded by attraction basins, while convergence to saddle points, if it is possible at all, is unstable with respect to small perturbations of the iterates. Another advantage of the MNM, as compared to the full Newton method, is that only $\frac14$ of the Hessian matrix has to be computed. In numerical experiments superior global convergence properties have been observed~\cite{Bakhurin2024}.

\medskip

In some cases, which are also relevant for training models in telecommunications, the mixed Hessian can be degenerate. This can be countered by adding a regularizing term, and it was proven in~\cite{Bakhurin2024} that when degeneration comes from a certain kind of symmetry group, as is the case in the mentioned application, a careful choice of this regularizer in theory does not alter the sequence of iterates. Another experimental observation was that far away from critical points the unregularized MNM can exhibit erratic behavior, but this can be alleviated by adding a regularization term.

In this paper, we consider the effects regularization has on the convergence properties of the MNM more thoroughly. We show that the statements of Theorem~\ref{thm:main} remain unaltered even with an arbitrary positive definite regularizing matrix. A second topic is the use of the MNM for the minimization of real analytic functions on $\mathbb R^n$. This can be done by considering them as function of complex variables after extending them to $\mathbb C^n$. Here we propose a regularization that fulfills at the same time the role of a penalty pushing the minimizers towards the real subspace $\mathbb R^n \subset \mathbb C^n$.

\section{Mixed Newton method with regularization}

Instead of the MNM iterate~\eqref{MNMiterate} we shall consider the iterate
\begin{equation} \label{RMNM}
z^{k+1} = z^k - \left( \frac{\partial^2f(z^k)}{\partial \bar z\partial z} + P \right)^{-1} \frac{\partial f(z^k)}{\partial\bar z},
\end{equation}
where $P \succ 0_{n\times n}$ is a fixed positive definite complex hermitian regularizing matrix. It assures that the regularized mixed Hessian is always positive definite and can also be used for damping the step size when chosen large enough. We shall now determine the behavior of~\eqref{RMNM} in the neighborhood of a critical point $\hat z$.

Recall that we minimize a function $f(z) = \sum_{j = 0}^{m - 1} |g_j(z)|^2$ with all $g_j$ holomorphic functions. Set $B = \frac{\partial^2f}{\partial\bar z\partial z} = \sum_{j = 0}^{m - 1} \overline{g_j'}(g_j')^\top$, $A = \sum_{j = 0}^{m - 1} g_j\overline{g_j''}$, where the derivatives are evaluated at the critical point $\hat z$ satisfying $\frac{\partial f}{\partial\bar z} = \sum_{j = 0}^{m - 1} g_j\overline{g_j'} = 0_n$.  Note that $B = B^* \succeq 0_{n\times n}$, $A = A^\top$, i.e., $B$ is complex-hermitian positive semi-definite, and $A$ is symmetric, but in general not hermitian.

Let us now recall the proof of Theorem~\ref{thm:main} in the case of the MNM iterate. Set $\delta = z - \hat z$, then we have for $z$ close to $\hat z$ that
\[ g_j(z) = g_j(\hat z) + (g_j'(\hat z))^\top\delta + O(\|\delta\|^2),
\]
\[ g_j'(z) = g_j'(\hat z) + g_j''(\hat z)\delta + O(\|\delta\|^2),
\]
and therefore
\begin{align*}
\frac{\partial f(z)}{\partial\bar z} &= \sum_{j = 0}^{m - 1} g_j(z)\overline{g_j'(z)} = \sum_{j = 0}^{m - 1} (g_j(\hat z) + (g_j'(\hat z))^\top\delta + O(\|\delta\|^2))(\overline{g_j'(\hat z) + g_j''(\hat z)\delta + O(\|\delta\|^2)}) \\
&= A\bar\delta + B\delta + O(\|\delta\|^2), \\
\frac{\partial^2f(z)}{\partial\bar z\partial z} &= \sum_{j = 0}^{m - 1} \overline{g_j'(z)}(g_j'(z))^\top = B + O(\|\delta\|).
\end{align*}

Then the regularized MNM iterate~\eqref{RMNM} acts on the differences $\delta_k = z^k - \hat z$ as
\begin{align*}
\delta_{k+1} &= \delta_k - \left( \frac{\partial^2f(z^k)}{\partial \bar z\partial z} + P \right)^{-1} \frac{\partial f(z^k)}{\partial\bar z} = \delta_k - \left( B + P + O(\|\delta\|) \right)^{-1} (A\bar\delta_k + B\delta_k + O(\|\delta_k\|^2)) \\
&= \delta_k - ( B + P )^{-1} (A\bar\delta_k + B\delta_k) + O(\|\delta_k\|^2) \\
&= \left( I - (B + P)^{-1}B \right) \delta_k - (B + P)^{-1}A \bar\delta_k + O(\|\delta_k\|^2).
\end{align*}
Setting $P = 0_{n\times n}$ we obtain the original MNM. Note that in this case, the first summand on the right-hand side disappears altogether if $B$ is invertible.

Let now $H_{\mathbb R} = \frac{\partial^2f(\hat z)}{\partial(\operatorname{Re}  (z), \operatorname{Im} (z))^2}$ be the Hessian at the critical point in real coordinates $(\operatorname{Re} (z), \operatorname{Im} (z))$. Then the signature of $H_{\mathbb R}$ is responsible for the type of the critical point. If $H_{\mathbb R} \succ 0_{2n\times 2n}$, then $\hat z$ is a minimum, and if the signature is mixed (with positive and negative eigenvalues), then the point $\hat z$ is a saddle point of the cost function $f$.

The matrix $H_{\mathbb R}$ can be transformed to the full Hessian matrix of Wirtinger derivatives by a conjugation with a non-degenerate coordinate change matrix. The Wirtinger derivatives, in turn, are given by the blocks $A, B$ defined above. We have
\[ M := \begin{pmatrix} \bar B & \bar A \\ A & B \end{pmatrix} = \begin{pmatrix} \frac{\partial^2f}{\partial z\partial\bar z} & \frac{\partial^2f}{\partial z^2} \\ \frac{\partial^2f}{\partial\bar z^2} & \frac{\partial^2f}{\partial \bar z\partial z} \end{pmatrix} = \frac14 \begin{pmatrix} I_n & -iI_n \\ I_n & iI_n \end{pmatrix} \begin{pmatrix} \frac{\partial^2f}{\partial \operatorname{Re} (z)^2} & \frac{\partial^2f}{\partial \operatorname{Re} (z) \, \partial \operatorname{Im} (z)} \\ \frac{\partial^2f}{\partial \operatorname{Im} (z) \, \partial \operatorname{Re} (z)} & \frac{\partial^2f}{\partial \operatorname{Im} (z)^2} \end{pmatrix} \begin{pmatrix} I_n & -iI_n \\ I_n & iI_n \end{pmatrix}^*,
\]
where $i^2 = -1$. It follows that the signature of the matrix $M$ on the left decides which type of critical point is represented by $\hat z$. Clearly, maxima are not possible, because by virtue of $B \succeq 0_{n\times n}$ the matrix $M$ cannot be negative definite.

We wish to put the signature of $M$ into relation with the stability of the dynamic system on the difference $\delta_k$. We have the extended dynamics
\[ \begin{pmatrix} \delta_{k+1} \\ \bar\delta_{k+1} \end{pmatrix} = \begin{pmatrix} I_n - (B + P)^{-1}B & - (B + P)^{-1}A \\ - (\bar B + \bar P)^{-1}\bar A & I_n - (\bar B + \bar P)^{-1}\bar B \end{pmatrix} \begin{pmatrix} \delta_k \\ \bar\delta_k \end{pmatrix} + O(\|\delta_k\|^2) =: Y\begin{pmatrix} \delta_k \\ \bar\delta_k \end{pmatrix} + O(\|\delta_k\|^2).
\]
Recall that the dynamics are stable (the iterates $\delta_k$ always converge to zero) if the spectrum of $Y$ is contained in the open unit disc.

The signature of a Hermitian matrix $M$ does not change under transformations of the type $M \mapsto CMC^*$, while the stability of the dynamical system $x_{k+1} = Yx_k$ defined by a matrix $Y$ does not change under transformations of the type $Y \mapsto FYF^{-1}$, where $C, F$ are arbitrary non-degenerate. Let $W = (B + P)^{1/2} = W^*$. Then $B = W^2 - P$. We have
\begin{align*}
\begin{pmatrix} \bar W^{-1} & 0_{n\times n} \\ 0_{n\times n} & W^{-1} \end{pmatrix} M \begin{pmatrix} \bar W^{-1} & 0_{n\times n} \\ 0_{n\times n} & W^{-1} \end{pmatrix} &= \begin{pmatrix} I_n - \bar W^{-1}\bar P\bar W^{-1} & \bar W^{-1}\bar AW^{-1} \\ W^{-1}A\bar W^{-1} & I_n - W^{-1}PW^{-1} \end{pmatrix} \\ &= \begin{pmatrix} I_n - \bar Q & \bar S \\ S & I_n - Q \end{pmatrix} =: M',
\end{align*}
\begin{align*}
\begin{pmatrix} 0_{n\times n} & \bar W \\ W & 0_{n\times n} \end{pmatrix} Y \begin{pmatrix} 0_{n\times n} & W^{-1} \\ \bar W^{-1} & 0_{n\times n} \end{pmatrix} &= \begin{pmatrix} \bar W^{-1}\bar P\bar W^{-1} & -\bar W^{-1}\bar AW^{-1} \\ -W^{-1}A\bar W^{-1} & W^{-1}PW^{-1} \end{pmatrix} \\ &= \begin{pmatrix} \bar Q & -\bar S \\ -S & Q \end{pmatrix} =: Y',
\end{align*}
where $S = W^{-1}A\bar W^{-1} = S^\top$, $Q = W^{-1}PW^{-1} = Q^* \succ 0_{n\times n}$. Note that the images of both $M, Y$ under these transforms are now Hermitian.

The dynamical system defined by $Y'$ is strictly stable if and only if its spectrum is contained in the interval $(-1,1)$. This is equivalent to the condition that the spectrum of $M' = I_{2n} - Y'$ is contained in the interval $(0,2)$, in which case $M' \succ 0_{2n\times 2n}$. Hence stability of the dynamical system defined by the RMNM in the vicinity of a critical point implies that the point is a minimum. This still leaves open the possibility that some minima are repelling, namely, if the spectrum of $M'$ is positive and its maximum eigenvalue exceeds 2.

Let us prove that this situation is not possible. Assume that $M' \succ 0_{2n\times 2n}$. Then the diagonal blocks of $M'$ are also positive definite. In particular, we get $0_{n\times n} \prec Q \prec I_n$, where the first inequality holds by the definition of $Q$. Hence $Q$ and $I_n-Q$ are contracting, and $(I_n-Q)^{-1/2},(I_n-\bar Q)^{-1/2}$ are expanding, i.e., application of these maps strictly increases the 2-norm of any non-zero vector. From $M' \succ 0_{2n\times 2n}$ it follows that $\sigma_{\max}((I_n-Q)^{-1/2}S(I_n-\bar Q)^{-1/2}) < 1$, i.e., this matrix product is contracting. But then $S$ must be contracting, and all 4 blocks of the matrix $M'$ are contracting. It follows that $\lambda_{\max}(M') = \sigma_{\max}(M') < 2$. The case of repelling minima is hence excluded.

We obtain the following result.

\begin{theorem}
The dynamical system defined by the regularized MNM in the neighborhood of a critical point is strictly stable if and only if this critical point is a minimum with positive definite Hessian. This dynamical system has repulsive directions if and only if the critical point is a saddle point. The case of a maximum is excluded.
\end{theorem}

\section{Minimization of real analytic functions}

In this section, we study how the MNM can be used to minimize a real analytic function $F$ on a simply connected domain $D_{\mathbb R} \subseteq \mathbb R^n$. The idea is to extend the function to a domain $D_{\mathbb C} \subseteq \mathbb C^n$ in complex space, to yield a function of the form $|g|^2$, where $g: D_{\mathbb C} \to \mathbb C$ is a holomorphic function.

Clearly, this is only possible if $F \geq 0$ on the domain $D_{\mathbb R}$. Since we may add constants to the objective, we may assume this condition without loss of generality. Then $F^2$ is also analytic on $D_{\mathbb R}$, moreover, minimizing $F^2$ is equivalent to minimizing $F$. If we define $g$ to be the holomorphic continuation of $F$ to a simply connected domain $D_{\mathbb C}$, then $F^2$ will be the restriction of a function $f: D_{\mathbb C} \to \mathbb R$ of the desired form $f = |g|^2$ with $g$ holomorphic.

This straightforward approach meets the following obstacles. The sum $\frac{\partial^2f}{\partial \bar z\partial z} = \sum_{j = 0}^{m - 1}\overline{\frac{dg_j}{dz}}\frac{dg_j}{dz}^\top$ defining the mixed Hessian consists of only one summand, hence the mixed Hessian is rank 1. Another problem is that the minima of $f$ can be located in the complex plane and a local minimum of $f$ on $D_{\mathbb R}$ may be a saddle point of $f$ on $D_{\mathbb C}$. The iterates will then be pushed away from $D_{\mathbb R}$.

Here we propose an approach that remedies both problems. It consists in adding to $f$ a sum of squares of holomorphic functions which penalizes deviations of $z$ into the complex plane but essentially does not alter the values of $f$ on $D_{\mathbb R}$. At the same time, the additional summands in the mixed Hessian regularize this matrix and render it positive definite.

Suppose that $g_0: D_{\mathbb C} \to \mathbb C$ a holomorphic function and that the goal is to minimize the modulus $|g_0|$ on the intersection $D_{\mathbb R} = \mathbb R^n \cap D_{\mathbb C}$. We shall minimize the function
\[ \sum_{l=0}^{2n} |g_l(z)|^2
\]
instead, where
\[ g_l(z) = \gamma e^{iz_l}, \quad l = 1,\dots,n; \qquad g_l(z) = \gamma e^{-iz_{(l - n)}}, \quad l = n+1,\dots,2n
\]
with $\gamma > 0$ a weighting coefficient. Then the additional term in the objective amounts to
\[ \sum_{l=1}^{2n} |g_l(z)|^2 = 2\gamma^2\sum_{l=1}^n \cosh(2\operatorname{Im}(z_l)).
\]
The additional term has the following properties:
\begin{itemize}
\item on $\mathbb R^n$, hence on $D_{\mathbb R}$ the constant $2n\gamma^2$ is added to the objective
\item for complex entries $z_l$ the modulus of the imaginary part is exponentially penalized
\item the amount of penalization can be tuned by choosing the coefficient $\gamma$
\end{itemize}

The mixed Hessian, which is initially of rank $1$, obtains an additive correction given by
\[ \sum_{l=1}^{2n} \overline{\frac{dg_l}{dz}}\frac{dg_l}{dz}^\top = 2\gamma^2\sum_{l=1}^n \cosh(2 \operatorname{Im}(z_l))e_le_l^\top \succeq 2\gamma^2I_n.
\]
Hence this correction acts as a positive definite regularizer.

Finally, we obtain the following iteration of the MNM for minimization of a nonnegative analytic function $F = g_0$ on $D_{\mathbb R}$ with additive complex-repulsive regularization:
\begin{equation}\label{RMNM_formula}
    z^{k+1} = z^k - \left( \overline{g_0'(z^k)}(g'_0(z^k))^\top + 2\gamma^2\sum_{j=1}^n \cosh(2 \operatorname{Im}(z^k_j))e_je_j^\top \right)^{-1} \left( g_0(z^k)\overline{g_0'(z^k)} + 2i\gamma^2\sum_{j=1}^n \sinh(2 \operatorname{Im}(z^k_j))e_j \right).
\end{equation}
Here $z^k_j$ is the $j$-th element of the current point $z^k$, and $e_j$ is the $j$-th basis vector.

Note that in the case when $g_0(\bar z) = \overline{g_0(z)}$, the Ordinary Newton Method (ONM) without regularization, acting in real space $\mathbb R^n$, proceeds by the iterate
\begin{equation}\label{ordinary_Newton_form}
    x^{k+1} = x^k - g_0(x^k) \left( g_0(x^k)g_0''(x^k) + g_0'(x^k)(g_0'(x^k))^\top\right)^{-1} g_0'(x^k).
\end{equation}
The RMNM, when started at a real point, will produce only real iterates in this case. Formula~\eqref{RMNM_formula} then simplifies to
\begin{equation}\label{RMNM_fromula_real}
    z^{k+1} = z^k - g_0(z^k)\left( g_0'(z^k)g'_0(z^k)^\top + 2\gamma^2I_n \right)^{-1}g_0'(z^k)  = z^k - \frac{g_0(z^k)g_0'(z^k)}{2\gamma^2 + \|g_0'(z^k)\|^2}.
\end{equation}

\section{Numerical experiments}

In this section, we present the results of our experiments for the comparison of the RMNM~\eqref{RMNM_formula} and the ONM~\eqref{ordinary_Newton_form}. We consider the problem of minimization of nonnegative real polynomials on $\mathbb{R}^2$.

\begin{example}\label{ex1_mins2}
Consider the polynomial
\begin{equation}\label{func_ex1_mins2}
  f(x) = f(x_1, x_2) = (2x_1 - 3x_2)^2  + x_1^2 (1-x_1)^2 + x_2^2 (1-x_2)^2.
\end{equation}

For this polynomial, we have a global minimum at the point $ (0,0)$ with optimal value $f^* = 0$. There is a local minimum at the point $\approx (1.04987, 0.709507)$ with value $ \approx 0.0460496$, and a saddle point at $ \approx (0.577876, 0.378919)$ with value $ \approx 0.11525$.
\end{example}

\begin{example}\label{ex2_mins2}
Consider the polynomial
\begin{equation}\label{func_ex2_mins2}
    f(x) = f(x_1, x_2) = (x_1 - x_2)^2 + x_1^2 (1-x_1)^2 + x_2^2 (2 - x_2)^2.
\end{equation}
For this polynomial, we have a global minimum at the point $ (0,0)$ with optimal value $f^* = 0$. There is a local minimum at the point $ \approx (1.27473, 1.81735)$ with value $ \approx 0.527264593$, and a saddle point at $(1,1 )$ with value $f(1,1) = 1$.
\end{example}

\begin{example}\label{ex3_mins3}
Consider the polynomial
\begin{equation}\label{func_ex3_mins3}
    f(x) = f(x_1, x_2) = (x_1 +1)^4 + (x_2 + 1)^4 + 4 x_1 x_2.
\end{equation}
For this polynomial, we have a global minimum at the point $ \approx (-0.31767219617, -0.31767219617)$ with optimal value $f^*  \approx 0.83717564078542$. There are local minima at the points $ \approx (0.1537213755, -1.53568738679)$, $ \approx 
 (-1.53568738679, 0.153721375541)$, both with value $ \approx 0.90983005625052$, and a saddle point at $(-1,0 )$.
\end{example}

Functions~\eqref{func_ex1_mins2}, \eqref{func_ex2_mins2} and~\eqref{func_ex3_mins3} are nonnegative and non-convex. For these functions, we run RMNM and ONM, with 625 different initial points in the square $[-1, 2] \times [-1,2]$ for Example~\ref{ex1_mins2}, with 1024 initial points in the square $[-1, 3] \times [-1, 3]$ for Example~\ref{ex2_mins2}, and with different 2601 initial points in the square $[-3, 2] \times [-3,2]$ for Example~\ref{ex3_mins3}. The initial points for both methods are the same and real. The results are presented in Table~\ref{table_1} and Table~\ref{table_2}, together with the used parameter value $\gamma$. These results show the number of initial points for which the methods converge to the local or global minimum and to the saddle point, they also show the number of initial points for which the methods diverge.

\begin{table}[htp]
\centering
\scalebox{0.72}{
\begin{tabular}{|c|cccc||cccc|}
\hline
\multirow{2}{*}{} & \multicolumn{4}{c||}{RMNM} & \multicolumn{4}{c|}{ONM}\\ \cline{2-9}
& \multicolumn{1}{c|}{ \shortstack{\# points  \\with converg. \\ to global min} } & \multicolumn{1}{c|}{\shortstack{\# points \\ with converg. \\ to local min}} & \multicolumn{1}{c|}{\shortstack{\# points \\ with converg. \\ to saddle }} & \shortstack{\# points \\ with \\ diverg.} & \multicolumn{1}{c|}{\shortstack{\# points  \\with converg. \\ to global min} } & \multicolumn{1}{c|}{\shortstack{\# points \\ with converg. \\to local min} } & \multicolumn{1}{c|}{\shortstack{ \# points \\ with converg. \\ to saddle}} & \shortstack{\# points \\ with \\ diverg. } \\ \hline
\shortstack{Example \ref{ex1_mins2}  \\with $625$ \\  initial points}& \multicolumn{1}{c|}{\textbf{625}} & \multicolumn{1}{c|}{0} & \multicolumn{1}{c|}{0} & 0 & \multicolumn{1}{c|}{276 } & \multicolumn{1}{c|}{319} & \multicolumn{1}{c|}{30} & 0  \\ \hline
\shortstack{Example \ref{ex2_mins2} \\ with $1024$ \\ initial points} & \multicolumn{1}{c|}{\textbf{1024}} & \multicolumn{1}{c|}{0} & \multicolumn{1}{c|}{0} &0  & \multicolumn{1}{c|}{463 } & \multicolumn{1}{c|}{443} & \multicolumn{1}{c|}{117} & 1 \\ \hline
\end{tabular}}
\caption{The results of RMNM with $\gamma = 10^{-3}$, and ONM for Examples~\ref{ex1_mins2} and~\ref{ex2_mins2}. }
\label{table_1}
\end{table}

\begin{table}[htp]
\centering
\begin{tabular}{|c|c|c|}
\hline
 &  RMNM &  ONM \\ \hline
\shortstack{\# points with convergence to the global min. } & \textbf{2601} & 672 \\ \hline
\shortstack{\# points with convergence to the first local min. } & 0 & 888 \\ \hline
\shortstack{\# points with convergence to the second local min. } & 0 & 888 \\ \hline
\shortstack{\# points with convergence to saddle point} & 0  & 76 \\ \hline
\shortstack{\# points without convergence \\ after $10^6$ iterations} &0  &77  \\ \hline
\end{tabular}
\caption{The results of RMNM with $\gamma = 10^{-2}$, and ONM for Example~\ref{ex3_mins3}, with 2601 initial points.}
\label{table_2}
\end{table}

From Tables~\ref{table_1} and~\ref{table_2}, we can see that RMNM always converges to the global minimum, with initialization by any initial point. This is due to the fact that the local minima become saddle points in complex space and are thus repelling for the method. In contrast to this, for ONM we observe convergence to the global and local minima and there are some initial points from which ONM converges to the saddle point or diverges, as for Example~\ref{ex2_mins2}.

\begin{remark}\label{remark_diverg}
If we take the initial points for RMNM of the form $(x_0 +  i, y_0 + i ) \in \mathbb{C}^2$ where $(x_0, y_0) \in [-1, 2] \times [-1, 2]$, (i.e., the imaginary parts of all coordinates are 1, not 0) then the RMNM will work worse. See Table~\ref{table_3}, it shows the results of the comparison between RMNM and ONM for Example~\ref{ex1_mins2}, where $(x_0, y_0) \in [-1, 2] \times [-1, 2] \subset \mathbb{R}^2$ are the initial points for ONM.

\begin{table}[ht]
\centering
\begin{tabular}{|c|c|c|}
\hline
 &  RMNM &  ONM \\ \hline
\shortstack{\# points with convergence to the global min. } & 25 & 24 \\ \hline
\shortstack{\# points with convergence to local min. } & 0 & 24 \\ \hline
\shortstack{\# points with convergence to saddle point} & 0  &1  \\ \hline
\shortstack{\# points without convergence \\ after $10^6$ iterations} &24  &0  \\ \hline
\end{tabular}
\caption{The results of RMNM with $\gamma = 10^{-3}$, and ONM for Example~\ref{ex1_mins2}, with 49 initial points, as in Remark~\ref{remark_diverg}.}
\label{table_3}
\end{table}

\end{remark}

\section{Conclusions}

In this work we modified the Mixed Newton Method and applied it to the optimization of real-valued functions of real variables as well as real-valued functions of complex variables. The modification consists of adding a regularization term to the mixed Hessian. This term, e.g., a multiple of the identity matrix in the case of a Levenberg-Marquardt type modification, alleviates degeneracies, e.g., caused by symmetries of the objective function. 

Alternatively, in the case of optimization of real parameters, regularization may be achieved by adding additional summands of a special kind in the objective function. These serve to push minima which lie in the complex space back to the real subspace, without essentially changing the restriction of the objective function to this real subspace.

The regularized mixed Newton method was tested on the task of minimizing non-convex real polynomials and showed superior global convergence properties by exhibiting an absolute preference for the global minimum. This is due to the fact that the local minima in real subspace turn to saddle points in complex space and thus become repelling for the method.

We compared the method on two training tasks for neural networks, one with application in telecommunications, the other a benchmark network. The regularized versions of MNM showed superior performance and required computational resources reduction in comparison with the classical competitors. These results are presented in the Appendix.

However, the regularized versions of MNM inherit from the MNM the limitation of its advantageous properties to the case of minimization of sums of squares of holomorphic (or analytic in the real case) functions.

\bibliography{neurips_2024}
\bibliographystyle{plain}

\newpage
\appendix
\section{Appendix}

\subsection{Applications}
At present, within the domain of machine learning applications, the prevailing methodology involves training models utilizing real-valued parameters on data that may be distributed in both complex and real spaces. In this respect, we introduce a novel approach characterized by transitioning from models with real-valued parameters to those with complex-valued parameters and training them using our proposed methods, no matter the distribution space of the data. Proposed methods based on MNM exhibit properties of repulsion from saddle point and reduced computational resources required for mixed Hessian calculation~(Theorem~\ref{thm:main}). In this context, the transition to models with complex-valued parameters, even when applied to real-valued data, ensures a high level of competitiveness compared to the approach of training models with real-valued parameters.

In this paper we consider two tasks to validate the proposed optimization methods and to demonstrate the advantages of transitioning to models with complex-valued parameters. The first one is a common problem in telecommunications, related to the non-linear distortion compensation generated by non-linear components in digital transceivers. Another problem is in essence a regression task, based on the common LIBSVM\footnote{\href{https://www.csie.ntu.edu.tw/~cjlin/libsvmtools/datasets/}{https://www.csie.ntu.edu.tw/~cjlin/libsvmtools/datasets/}} dataset. 

Both tasks are considered in real and complex space, i.e., for both issues two models are provided, one with real and one with complex parameters. 

The complex models are trained by 4 algorithms: mixed Newton method and Newton method with Levenberg-Marquardt type adaptive regularization control (LM-MNM and LM-NM), cubic Newton method (CNM) and cubic mixed Newton method (CMNM)~\cite{Levenberg1944, Marquardt1963, Nesterov2006}. At the same time, the real-valued models are optimized by 2 algorithms: the Levenberg-Marquardt regularized Newton method (LM-NM) and cubic Newton method (CNM).

In addition, we assume 3 different starting points for the model with complex parameters: sampled from the real axis, the imaginary axis, and the whole complex plane.

\begin{figure}[ht]
  \centering
  \begin{minipage}[b]{0.48\textwidth}
    \centering
    \includegraphics[width=\textwidth]{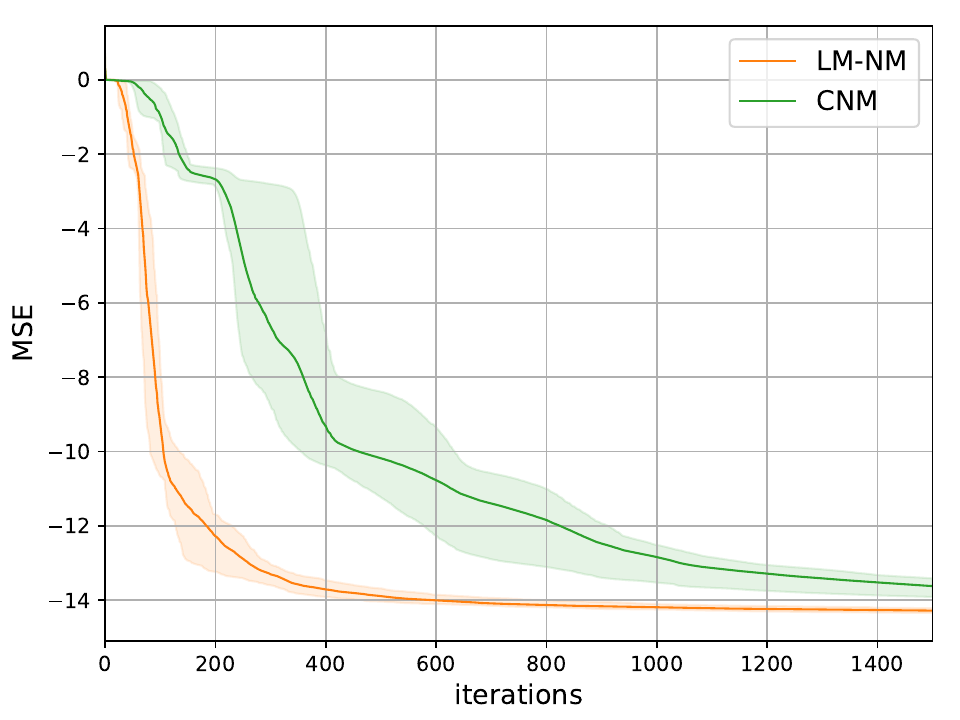}
  \end{minipage}
  \begin{minipage}[b]{0.48\textwidth}
    \centering
    \includegraphics[width=\textwidth]{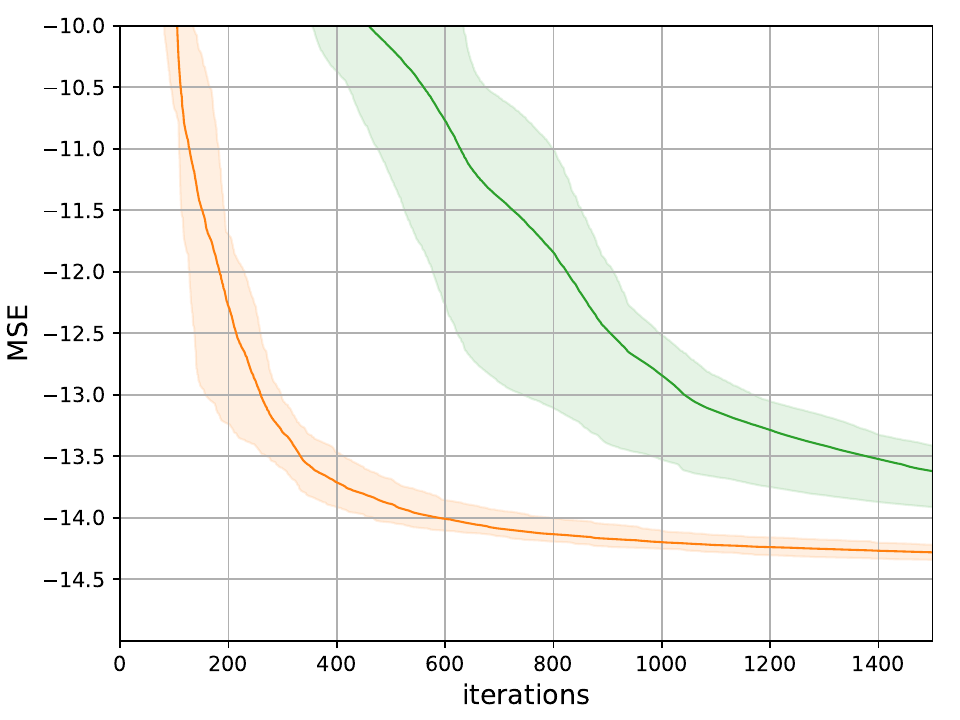}
  \end{minipage}
  \captionsetup{justification=centering}
  \caption{Learning curves for the tested optimization algorithms. RV-CNN model.}
  \label{fig:rvcnn}
\end{figure}

Note, that in Figures~\ref{fig:cvcnn_complex}--~\ref{fig:rvcnn} and~\ref{fig:abalone}
each learning curve corresponds to the average value of trials with five initial points and the fill reflects the min-max range:
\begin{equation}
    \displaystyle\text{MSE}_{\text{aver}}=\frac{1}{N}\sum_{j=0}^{N-1}\text{MSE}_{j}, \quad
    \text{MSE}_{j}\in\mathbb{R}^{T}, \nonumber
\end{equation}
\begin{align}
    \displaystyle\text{MSE}_{\text{min}}=\min_{j}(\text{MSE}_{j}), \quad \displaystyle\text{MSE}_{\text{max}}=\max_{j}(\text{MSE}_{j}), \quad j=\overline{0, N-1},
    \label{mse_agregated}
\end{align}
where number of trials $N=5$, number of iterations $T$, differ for telecommunications and LIBSVM task. Normalized MSE in~Figures~\ref{fig:cvcnn_complex}--~\ref{fig:compare} is calculated on base of~\eqref{mse_agregated}. For instance $\text{NMSE}_{\text{aver}}~=~20\log_{10}\text{MSE}_{\text{aver}}$.

At the same time, in Figures~\ref{fig:compare}~and~\ref{fig:compare_abalone} learning curves for model with complex weights are aggregated among all initial points trials:
\begin{equation}
    \displaystyle\text{MSE}_{\text{aver}}=\frac{1}{P\cdot N}\sum_{p=0}^{P-1}\sum_{j=0}^{N-1}\text{MSE}_{j, p}, \quad
    \text{MSE}_{j, p}\in\mathbb{R}^{T}, \nonumber
\end{equation}
\begin{align}
    \displaystyle\text{MSE}_{\text{min}}=\min_{j, p}(\text{MSE}_{j, p}), \quad \displaystyle\text{MSE}_{\text{max}}=\max_{j, p}(\text{MSE}_{j, p}), \quad j=\overline{0, N-1}, \quad p=\overline{0, P-1},
    \label{mse_agregated_start_points}
\end{align}
where $P=3$ -- corresponds to 3 different initial points sampling cases: from real axis, imaginary axis and whole complex plane.

\subsection{Derivation of the iterations of LM-MNM and CMNM}

The Newton method with Levenberg-Marquardt adaptive regularization and the cubic Newton method are well-known~\cite{Levenberg1944, Marquardt1963, Nesterov2006}. Here we derive a complex version of these methods to modify the MNM.

Introduce variables
\[ z = x + iy = \operatorname{Re}(z) + i\operatorname{Im}(z),\quad \bar{z} = x - iy,\quad r = (x, y) = (\operatorname{Re}(z), \operatorname{Im}(z)),\quad c = (z, \bar{z}).
\]
Then the gradients and Hessians with respect to $r$ and $c$ are linearly dependent on each other, as these variables can be transformed one into another by a linear transformation. Denote $\mathcal{H}_{\bar{z}z} = \frac{\partial}{\partial \bar{z}} \left(\frac{\partial f}{\partial z}\right)^{\top}$ etc. Then we have
$$
\mathcal{H}_{rr} = J^{H} \mathcal{H}_{cc} \overline{J} \Longrightarrow \mathcal{H}_{cc} = \frac{1}{4} J\mathcal{H}_{rr} J^{\top},
$$
$$\text{where }J = \begin{pmatrix}
 I_n& -iI_n \\
I_n & iI_n
\end{pmatrix}.$$
Consider the second-order Taylor approximation of the real-valued target function $f$ of complex variables $z = x + iy$ as a function of $r = (x, y)$ near a real-valued point $r$:
\begin{equation*}
    \begin{aligned}
        &M_{r}^{2}(\Delta_{r}) := f(r) + \left\langle g_{r}, \Delta_{r}\right\rangle + \frac{1}{2}\left\langle\mathcal{H}_{rr}\Delta_{r}, \Delta_{r}\right\rangle\in\mathbb{R},\quad\Delta_{r} := r_{1} - r\in\mathbb{R}^{2n},\quad g_r := \frac{\partial f}{\partial r}.
    \end{aligned}
\end{equation*}
We impose the following Lipschitz property on the target function's gradient and Hessian:
\begin{equation}\label{eq:lip}
    \begin{aligned}
        &\forall r_{1}, r\in\mathbb{R}^{2n},~\exists L_{g} > 0:\quad\|g_{r_{1}} - g_{r}\|\leq L_{g}\|r_{1} - r\|;\\
        &\forall r_{1}, r\in\mathbb{R}^{2n},~\exists L_{\mathcal{H}} > 0:\quad\|\mathcal{H}_{r_{1}r_{1}} - \mathcal{H}_{rr}\|\leq L_{\mathcal{H}}\|r_{1} - r\|.
    \end{aligned}
\end{equation}
Define $\Delta_{c} = c_{1} - c\in\mathbb{C}^{2n}$ and $\Delta_z = z_1 - z\in\mathbb{C}^n$, then $\|\Delta_{r}\|^{2} = \frac{1}{2}\|\Delta_{c}\|^{2}$. The Taylor polynomial can be rewritten as
\begin{equation*}
    \begin{aligned}
        &M_{c}^{2}\left(\Delta_{c}\right) := f(c) + \left\langle g_{\bar{c}}, \Delta_{c}\right\rangle + \frac{1}{2}\left\langle\mathcal{H}_{\bar{c}c}\Delta_{c}, \Delta_{c}\right\rangle\in\mathbb{R}.
    \end{aligned}
\end{equation*}
To build the correct procedure with line search of the Lipschitz constant, one should consider the original bound to check $L_{\mathcal{H}}\geq0$:
\begin{equation*}
    \begin{aligned}
        f(z_1)&\leq f(z) + \left\langle g_{\bar{c}}, \Delta_c\right\rangle + \frac{1}{2}\left\langle \mathcal{H}_{\bar{c}c}\Delta_c, \Delta_c\right\rangle + \frac{L_{\mathcal{H}}}{12\sqrt{2}}\left\|\Delta_c\right\|^3 = \\
        &= f(z) + 2\left\langle g_{\bar{z}}, \Delta_{z}\right\rangle + \left\langle\mathcal{H}_{\bar{z}z}\Delta_{z}, \Delta_{z}\right\rangle + \operatorname{Re}\left(\left\langle \mathcal{H}_{zz}\Delta_{z}, \Delta_{z}\right\rangle\right) + \frac{L_{\mathcal{H}}}{6}\left\|\Delta_{z}\right\|^{3} = \\
        &= f(z) + 2\left\langle g_{\bar{z}}, \Delta_{z}\right\rangle + \left\langle\mathcal{H}_{\bar{z}z}\Delta_{z}, \Delta_{z}\right\rangle + \frac{1}{4}\left\langle \left(\mathcal{H}_{xx} - \mathcal{H}_{yy}\right)\Delta_{z}, \Delta_{z}\right\rangle + \frac{L_{\mathcal{H}}}{6}\left\|\Delta_{z}\right\|^{3} = \\
        &= f(z) + 2\left\langle g_{\bar{z}}, \Delta_{z}\right\rangle + \left\langle\mathcal{H}_{\bar{z}z}\Delta_{z}, \Delta_{z}\right\rangle + \frac{1}{2}\left\langle \left(\mathcal{H}_{zz} + \mathcal{H}_{\bar{z}\bar{z}}\right)\Delta_{z}, \Delta_{z}\right\rangle + \frac{L_{\mathcal{H}}}{6}\left\|\Delta_{z}\right\|^{3} = \\
        &= f(z) + \left\langle \left(2\cdot g_{\bar{z}}\right), \Delta_{z}\right\rangle + \frac{1}{2}\cdot\left\langle\left(2\cdot \mathcal{H}_{\bar{z}z} + \mathcal{H}_{zz} + \mathcal{H}_{\bar{z}\bar{z}}\right)\Delta_{z}, \Delta_{z}\right\rangle + \frac{L_{\mathcal{H}}}{6}\cdot\left\|\Delta_{z}\right\|^{3}.
    \end{aligned}
\end{equation*}
Here we used that
\begin{equation*}
    \begin{aligned}
        &\operatorname{Re}\left(\left\langle\mathcal{H}_{zz}\Delta_z, \Delta_z\right\rangle\right) = \frac{1}{2}\left(\left\langle\mathcal{H}_{zz}\Delta_z, \Delta_z\right\rangle + \left\langle\mathcal{H}_{\bar{z}\bar{z}}\Delta_{\bar{z}}, \Delta_{\bar{z}}\right\rangle\right) = \frac{1}{4}\left\langle\left(\mathcal{H}_{xx} - \mathcal{H}_{yy}\right)\Delta_z, \Delta_z\right\rangle = \frac{1}{2}\left\langle\left(\mathcal{H}_{zz} + \mathcal{H}_{\bar{z}\bar{z}}\right)\Delta_z, \Delta_z\right\rangle.
    \end{aligned}
\end{equation*}
Finally, the correct Cubic Newton Method step is validated by the following conditions:
\begin{equation*}
    \begin{aligned}
        &f(z_1)\leq \underbrace{f(z) + 2\left\langle g_{\bar{z}}, \Delta_{z}\right\rangle + \left\langle\left(\mathcal{H}_{\bar{z}z} + \frac{1}{2}\left(\mathcal{H}_{zz} + \mathcal{H}_{\bar{z}\bar{z}}\right)\right)\Delta_{z}, \Delta_{z}\right\rangle + \frac{L_{\mathcal{H}}}{6}\delta^{3}}_{L_{\mathcal{H}}\text{ is selected first, bound is computed after finding }\delta\text{ and computing }z_1};\\
        &\delta^2 = \left\|\left(\mathcal{H}_{\bar{z}z} + \frac{1}{2}\left(\mathcal{H}_{zz} + \mathcal{H}_{\bar{z}\bar{z}}\right) + \frac{\delta L_{\mathcal{H}}}{4}I_{n}\right)^{-1}g_{\bar{z}}\right\|^2,~\delta\geq\frac{4}{L_{\mathcal{H}}}\max\left\{-\lambda_{\min}\left(\mathcal{H}_{\bar{z}z} + \frac{1}{2}\left(\mathcal{H}_{zz} + \mathcal{H}_{\bar{z}\bar{z}}\right)\right), 0\right\};\\
        &z_1 := z - \left(\mathcal{H}_{\bar{z}z} + \frac{1}{2}\left(\mathcal{H}_{zz} + \mathcal{H}_{\bar{z}\bar{z}}\right) + \frac{\delta L_{\mathcal{H}}}{4}I_{n}\right)^{-1}g_{\bar{z}};\\
        &c_1 := \begin{pmatrix}
            z_1\\
            \bar{z}_1
        \end{pmatrix}.
    \end{aligned}
\end{equation*}
The complex Hessian in the Cubic Newton Method update rule can be easily expressed via a linear combination of real-valued Hessians:
\begin{equation*}
    \begin{aligned}
        \mathcal{H}_{\bar{z}z} + \frac{1}{2}\left(\mathcal{H}_{zz} + \mathcal{H}_{\bar{z}\bar{z}}\right) + \frac{\delta L_{\mathcal{H}}}{4}I_n &= \frac{1}{4}\left(\mathcal{H}_{xx} + \mathcal{H}_{yy} + i\left(\mathcal{H}_{yx} - \mathcal{H}_{xy}\right)\right) + \frac{1}{2}\left( \frac{1}{4}\left(\mathcal{H}_{xx} - \mathcal{H}_{yy} - i\left(\mathcal{H}_{xy} + \mathcal{H}_{yx}\right)\right) + \right.\\
        &\left. + \frac{1}{4}\left(\mathcal{H}_{xx} - \mathcal{H}_{yy} + i\left(\mathcal{H}_{yx} + \mathcal{H}_{xy}\right)\right)\right)  + \frac{\delta L_{\mathcal{H}}}{4}I_n =\\
        &= \frac{1}{4}\cdot\left(\left(2\cdot\mathcal{H}_{xx} + \delta\cdot L_{\mathcal{H}}\cdot I_n\right) + i\cdot\left(\mathcal{H}_{yx} - \mathcal{H}_{xy}\right)\right);\\
        z_1 &:= z - 4\cdot\left(2\cdot\mathcal{H}_{xx} + \delta\cdot L_{\mathcal{H}}\cdot I_n + i\cdot\left(\mathcal{H}_{yx} - \mathcal{H}_{xy}\right)\right)^{-1}g_{\bar{z}} =\\
        &= z - 2\cdot\left(2\cdot\mathcal{H}_{xx} + \delta\cdot L_{\mathcal{H}}\cdot I_n + i\cdot\left(\mathcal{H}_{yx} - \mathcal{H}_{xy}\right)\right)^{-1}\left(g_x + i\cdot g_y\right);\\
        c_1 &:= \begin{pmatrix}
            z_1\\
            \bar{z}_1
        \end{pmatrix}.
    \end{aligned}
\end{equation*}
Note that there is no need in $\mathcal{H}_{yy}$ computation. To construct the main step of the Cubic Mixed Newton Method (CMNM) we use the following inequality:
\begin{equation*}
    \begin{aligned}
        f(z_{1})&\leq f(z) + \left\langle \left(2\cdot g_{\bar{z}}\right), \Delta_{z}\right\rangle + \frac{1}{2}\cdot\left\langle\left(2\cdot \mathcal{H}_{\bar{z}z}\right)\Delta_{z}, \Delta_{z}\right\rangle + \left\langle\frac{1}{2}\cdot\left(\mathcal{H}_{zz} + \mathcal{H}_{\bar{z}\bar{z}}\right)\Delta_{z}, \Delta_{z}\right\rangle + \frac{L_{\mathcal{H}}}{6}\cdot\left\|\Delta_{z}\right\|^{3}\leq\\
        &\leq f(z) + \left\langle \left(2\cdot g_{\bar{z}}\right), \Delta_{z}\right\rangle + \frac{1}{2}\left\langle\left(2\cdot\mathcal{H}_{\bar{z}z}\right)\Delta_{z}, \Delta_{z}\right\rangle + L_g\cdot\left\|\Delta_z\right\|^2 + \frac{L_{\mathcal{H}}}{6}\cdot\left\|\Delta_{z}\right\|^{3}\leq\\
        &\leq f(z) + \left\langle \left(2\cdot g_{\bar{z}}\right), \Delta_{z}\right\rangle + \frac{1}{2}\left\langle\left(2\cdot\mathcal{H}_{\bar{z}z}\right)\Delta_{z}, \Delta_{z}\right\rangle + \hat{L}\cdot\left\|\Delta_z\right\|^2 + \frac{\hat{L}}{6}\cdot\left\|\Delta_{z}\right\|^{3},\quad \hat{L}\geq\max\left\{L_{g}, L_{\mathcal{H}}\right\}.
    \end{aligned}
\end{equation*}
With this bound one may construct the CMNM update rule:
\begin{equation*}
    \begin{aligned}
        &f(z_{1})\leq f(z) + 2\left\langle g_{\bar{z}}, \Delta_{z}\right\rangle + \left\langle\mathcal{H}_{\bar{z}z}\Delta_{z}, \Delta_{z}\right\rangle + \hat{L}\delta^{2}\left(1 + \frac{\delta}{6}\right);\\
        &\delta^{2} = \left\|\left(\mathcal{H}_{\bar{z}z} + \hat{L}\left(1 + \frac{\delta}{4}\right)I_{n}\right)^{-1}g_{\bar{z}}\right\|^{2},~\delta\geq4\max\left\{\frac{-\lambda_{\min}\left(\mathcal{H}_{\bar{z}z}\right)}{\hat{L}} - 1,~0\right\};\\
        &z_{1} := z - \left(\mathcal{H}_{\bar{z}z} + \hat{L}\left(1 + \frac{\delta}{4}\right)I_{n}\right)^{-1}g_{\bar{z}}.
    \end{aligned}
\end{equation*}

In the simulations the algorithms MNM and NM are accelerated by means of Levenberg-Marquardt adaptive regularization control. Let us denote by $\mathcal{H}$ the Hessian matrix calculated for the appropriate algorithm: MNM or NM. Then the adaptive $\lambda$-regularization control is formalized in Algorithm~\ref{alg:lm_adaptive}.

\begin{algorithm}[ht]
\caption{Levenberg-Marquardt adaptive regularization control.}\label{alg:lm_adaptive}
 \KwData{$x$ -- \text{input signal (input data)}, $y$ -- \text{desired output signal (target variable(s), label(s))}, \textbf{model}, $f$ -- loss evaluation\;}
 \KwResult{$z$ -- model parameters\;}
 Initialize regularization parameter $\lambda\geq0$\;
 Initialize hyper-parameters $\alpha > 1, \mu > 0$\;
 Evaluate initial loss: $f:=\text{loss}(\textbf{model}(x, z), y)$\;
 \While{Training stop criteria not achieved}{
 Calculate the appropriate Hessian $\mathcal{H}$ and the appropriate gradient $\nabla(f)$\;
 Store current model parameters: $z_{tmp}:=z$\;
 Recalculate model parameters $z:= z - \mu (\mathcal{H} + \lambda\|\operatorname{vec}(\mathcal{H})\|_{\infty} I)^{-1}\nabla(f)$\;
 Calculate optimization criterion: $f_{1}:=\text{loss}(\textbf{model}(x, z), y)$\;
\eIf{$f_{1} < f$}{
$\lambda := \frac{\lambda}{\alpha}$\;
}{
\While{$f_{1} \geq f$}{
  $\lambda := \lambda \cdot \alpha$\;
  $z:= z_{tmp} - \mu (\mathcal{H} + \lambda\|\operatorname{vec}(\mathcal{H})\|_{\infty} I)^{-1}\nabla(f)$\;
  $f_{1}:=\text{loss}(\textbf{model}(x, z), y)$\;
}
 } $f := f_1$\;
 }
\end{algorithm}

\subsubsection{Digital Pre-Distortion}
Complex-valued neural networks are used in nonlinear signal processing of wireless communication systems while processing baseband signals, which are presented in complex form~\cite{Luo2011}. One of the most important tasks in nonlinear digital processing is the digital predistorting for transmitter power amplifier (PA) linearization. At the same time, the development of wireless systems requires simultaneous processing of several transmission frequency bands, as a result of which there is a need for multi-band nonlinear processing. For two bands we have dualband predistorting~\cite{Ghannouchi2015, Molina2017}.

The baseline solution for this task is a memory polynomial approximation of the PA output $y_{1_n}$ and $y_{2_n}$ as a function of the input signals $x_{1_n}$ and $x_{2_n}$ in the parameterized form
\begin{equation}
y_{1_n} = x_{1_n} - \sum_{d = -D}^D\sum_{m=0}^P\sum_{k=0}^P z_{1_{dmk}} x_{1_{n-d}} |x_{1_{n-d}}|^m|x_{2_{n-d}}|^k;
\end{equation}
\begin{equation}
y_{2_n} = x_{2_n} - \sum_{d = -D}^D\sum_{m=0}^P\sum_{k=0}^P z_{2_{dmk}} x_{2_{n-d}} |x_{2_{n-d}}|^m|x_{1_{n-d}}|^k,
\end{equation}
where $n$ is the discrete time, $d$ are time delays, $z_{1_{dmk}}$ and $z_{2_{dmk}}$ are complex coefficients which must be estimated for the first and second bands by training. The main baseline model is a single-layer structure whose coefficients can be optimized by the least squares (LS) approach. The baseline models main disadvantage is the quadratic parameter number growth due to exploitation of 2-dimensional nonlinear functions. The total number of coefficients of the baseline model is $2(D+1)(P+1)^2$.  

In contrast, complex multi-layer neural networks are found to have less complexity while maintaining high performance. In the present work we consider a Complex-Valued Convolutional Neural Network (CV-CNN) and a Real-Valued Convolutional Neural Network (RV-CNN).

These models are fed by the following matrices, correspondingly:
\begin{equation}
    \wave{x}_n = 
        \displaystyle\begin{pmatrix}
            x_n, \\
            |x_n|
        \end{pmatrix}\in\mathbb{C}^{2\times L}\, \text{  for CV-CNN}, \, \text{  
        }\wave{x}_n =
        \displaystyle\begin{pmatrix}
            |x_n|\\
            \text{Re}(x_n)\\
            \text{Im}(x_n)
        \end{pmatrix}\in\mathbb{R}^{3\times L}\text{  for RV-CNN},
\end{equation}
where $L$ is the length of the input signal $x_n\in\mathbb{C}^{1\times L}$ at time instant $n$, here operation $|\cdot|$ is applied element-wise.

Denote by $\mathbb{K}$ the real field $\mathbb{R}$ or the complex field $\mathbb{C}$, depending on the considered model. Then, the output of the $p$-th convolutional layer of the CV-CNN or RV-CNN can be expressed as follows:
\begin{align}
    u_{n,p}=
    \begin{pmatrix}\displaystyle\sigma(b_{0,p}+\sum_{c=0}^{K-1}u_{n,p-1,c} * z_{c,0,p}) \\
    \displaystyle\sigma(b_{1,p}+\sum_{c=0}^{K-1}u_{n,p-1,c} * z_{c,1,p}) \\
    \vdots \\
    \displaystyle\sigma(b_{K-1,p}+\sum_{c=0}^{K-1}u_{n,p-1,c} * z_{c,K-1,p})
     \end{pmatrix}\in\mathbb{K}^{K\times L},
\end{align}
where $\mathbb{K}=\mathbb{R}$ in case of RV-CNN, $\mathbb{K}=\mathbb{C}$ in case of CV-CNN, $\sigma(\cdot)$ -- sigmoid function. In addition, $c$ is the input channel index, $p$ the layer index, $n$ the time index, $z_{c,k,p}\in\mathbb{K}^{1\times M}$ the 1D convolution kernel of the input channel $c$, the output channel $k$, the layer $p$. Further $M$ is the kernel size, $b_{c,p}\in\mathbb{K}$ the bias, $u_{n,p-1,c}\in\mathbb{K}^{1\times L}$ the $c$-th output of convolutional layer $p-1$. Note that the input of the $0$-th convolutional layer is formally given by the model input channel, $u_{n,-1,c}=\wave{x}_{n,c}\in\mathbb{K}^{1\times L}$.

Thus, the output of the whole RV-CNN can be described as
\begin{align}
    \begin{pmatrix}
        u_{n, Re} \\ u_{n, Im}
    \end{pmatrix}=\sigma\begin{pmatrix}
    \displaystyle b_{0,\text{out}}+\sum_{c=0}^{K-1}u_{n,P-1,c} * z_{c,0,\text{out}} \\
    \displaystyle b_{1,\text{out}}+\sum_{c=0}^{K-1}u_{n,P-1,c} * z_{c,1,\text{out}}
  \end{pmatrix}\in\mathbb{R}^{2\times L},
\end{align}
$z_{c,0,\text{out}}$, $z_{c,1,\text{out}}\in\mathbb{R}^{1\times M}$, $b_{0,\text{out}}$, $b_{1,\text{out}}\in\mathbb{R}$. Finally, the complex output vector is formed as
\begin{align}
    y_n=u_{n, Re} + iu_{n, Im}\in\mathbb{C}^{1\times L}.
\end{align}

For the CV-CNN model the output is given by
\begin{align}
    y_{n}=\sigma\begin{pmatrix}
    \displaystyle b_{\text{out}}+\sum_{c=0}^{K-1}u_{n,P-1,c} * z_{c,\text{out}}
  \end{pmatrix}\in\mathbb{C}^{1\times L},
\end{align}
where $z_{c,\text{out}}\in\mathbb{C}^{1\times M}$, $b_{\text{out}}\in\mathbb{C}$. 

In current simulations we assume hyper-parameters of the models which are presented in Table~\ref{tab:table_telecom_hyper}. The train signal block length $x_n\in\mathbb{C}^{1\times L}$ in the simulations is set to $L=170800$.
\begin{table}[ht]
\centering
\caption{Hyper-parameters of the models in the telecommunications application.}
\begin{tabular}{|c|c|c|c|c|c|}
\hline
 & \shortstack{kernel \\ size, $M$} & 
 \shortstack{layer \\ num., $P$} & \shortstack{output \\ num., $K$} &
 \shortstack{space, \\ $\mathbb{K}$} &
 \shortstack{param., \\ num.} \\ \hline
\shortstack{RV-CNN} & 3 & 4 & 6, 5, 5, 2 & $\mathbb{R}$ & 321 \\ \hline
\shortstack{CV-CNN} & 3 & 4 & 3, 3, 3, 1 & $\mathbb{C}$ & 109 \\ \hline
\end{tabular}
\label{tab:table_telecom_hyper}
\end{table}

\begin{figure}[ht]
  \centering
  \begin{minipage}[b]{0.48\textwidth}
    \centering
    \includegraphics[width=\textwidth]{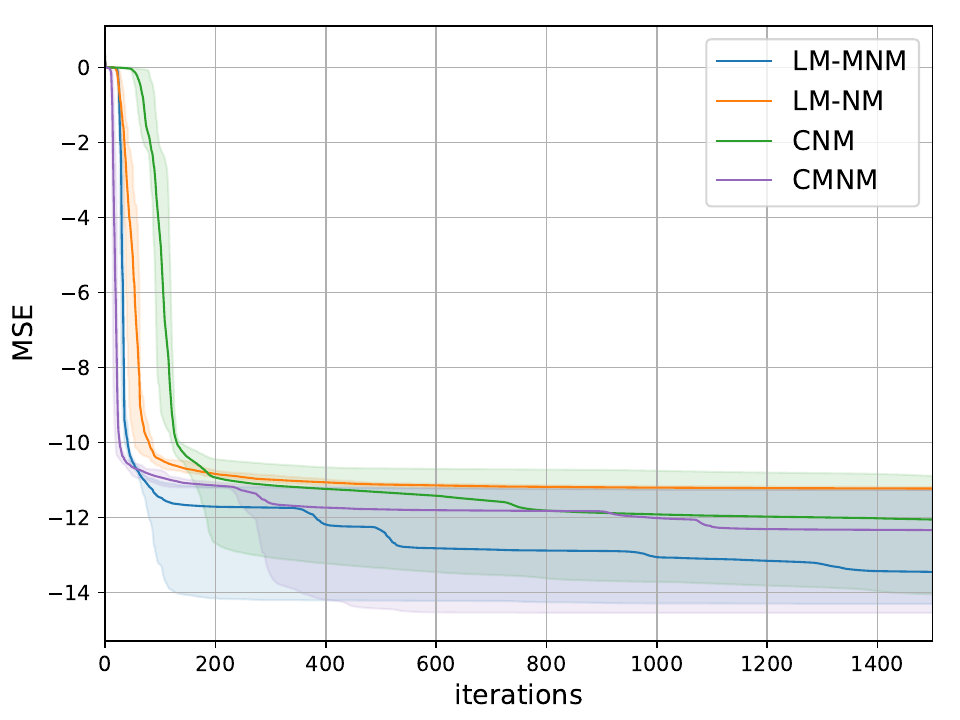}
  \end{minipage}
  \begin{minipage}[b]{0.48\textwidth}
    \centering
    \includegraphics[width=\textwidth]{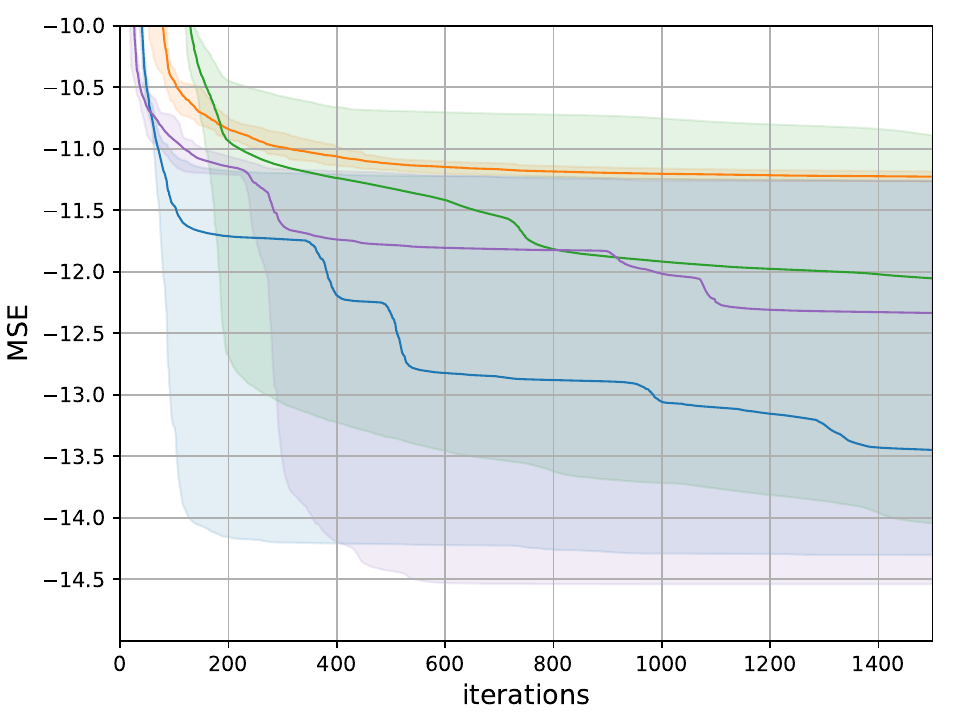}
  \end{minipage}
  \captionsetup{justification=centering}
  \caption{Learning curves for the tested optimization algorithms. CV-CNN model. Initial parameters in the complex plane.}
  \label{fig:cvcnn_complex}
\end{figure}

\begin{figure}[ht]
  \centering
  \begin{minipage}[b]{0.48\textwidth}
    \centering
    \includegraphics[width=\textwidth]{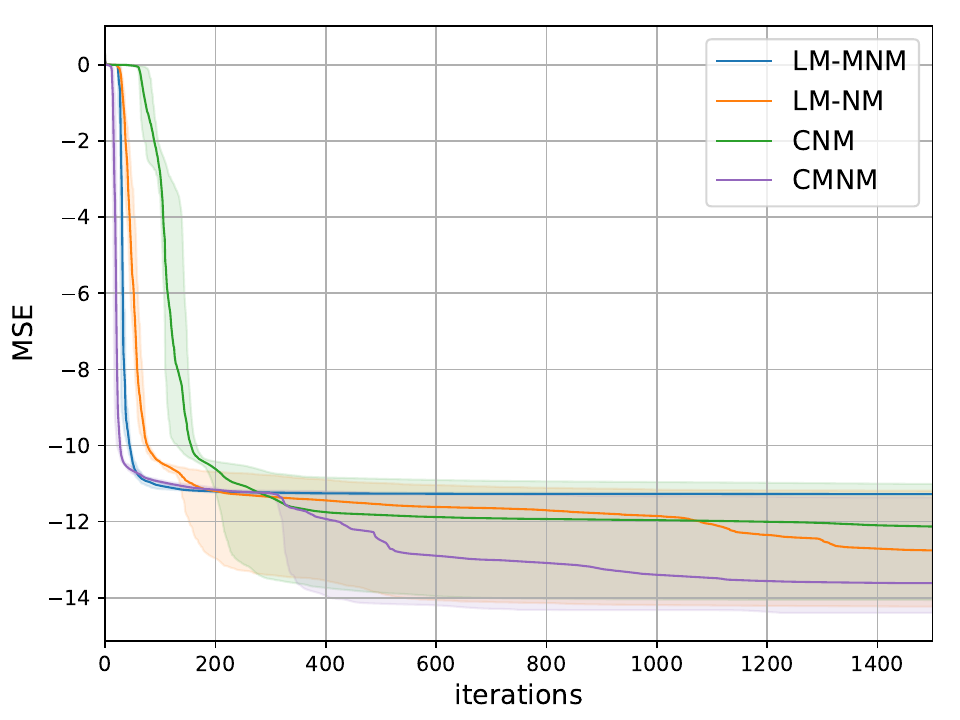}
  \end{minipage}
  \begin{minipage}[b]{0.48\textwidth}
    \centering
    \includegraphics[width=\textwidth]{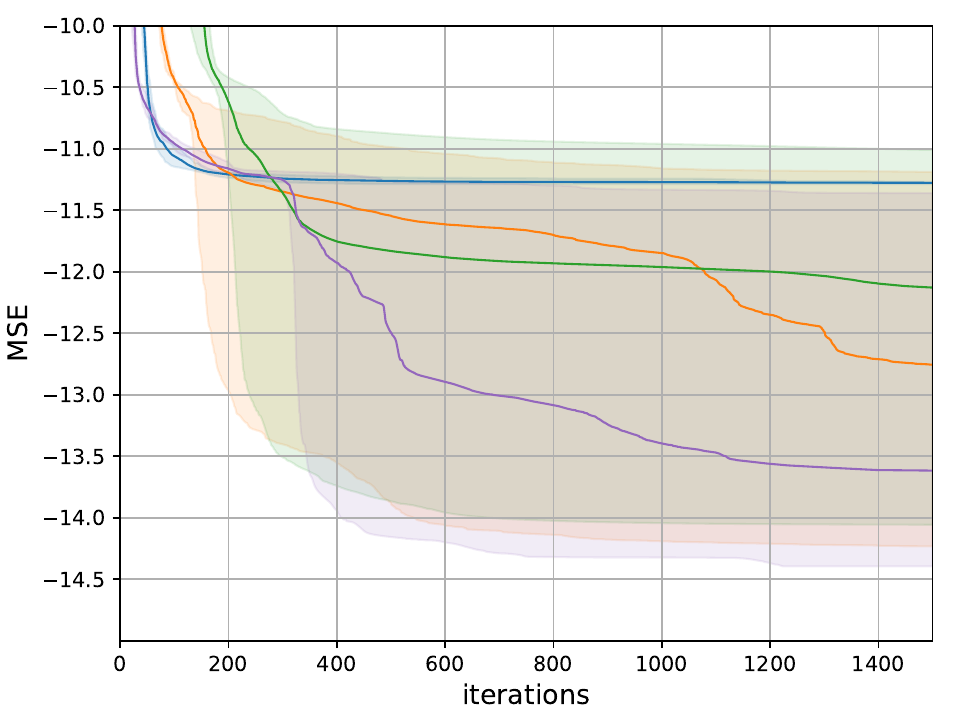}
  \end{minipage}
  \captionsetup{justification=centering}
  \caption{Learning curves for the tested optimization algorithms. CV-CNN model. Initial parameters on the real axis.}
  \label{fig:cvcnn_real}
\end{figure}

\begin{figure}[ht]
  \centering
  \begin{minipage}[b]{0.48\textwidth}
    \centering
    \includegraphics[width=\textwidth]{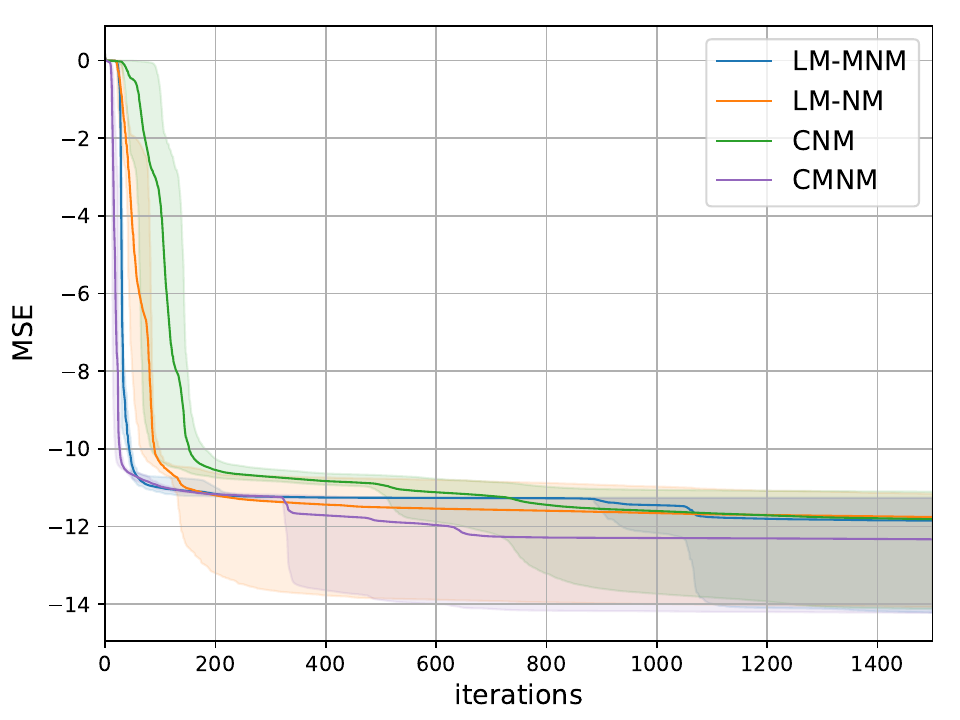}
  \end{minipage}
  \begin{minipage}[b]{0.48\textwidth}
    \centering
    \includegraphics[width=\textwidth]{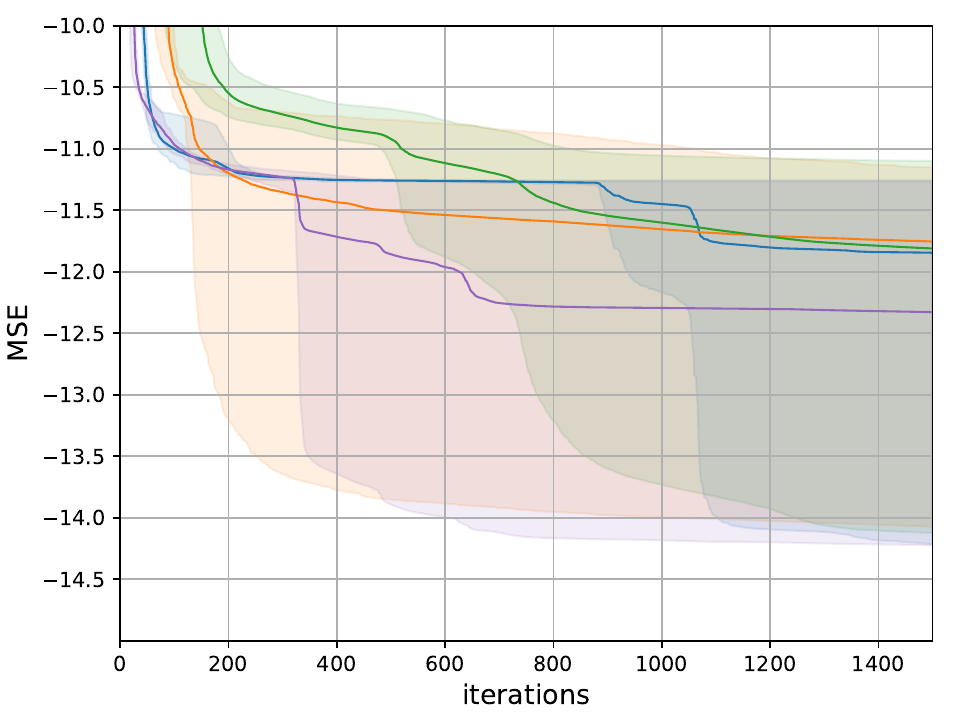}
  \end{minipage}
  \captionsetup{justification=centering}
  \caption{Learning curves for the tested optimization algorithms. CV-CNN model. Initial parameters on the imaginary axis.}
  \label{fig:cvcnn_imag}
\end{figure}

In Figures~\ref{fig:cvcnn_complex}, \ref{fig:cvcnn_real} and~\ref{fig:cvcnn_imag} the four proposed second-order methods are compared in terms of convergence speed and performance. It can be seen that the LM-MNM and CMNM yield the best normalized mean square error (NMSE) values after 1500 iterations for each starting point, except for the case of initial parameters sampled on the real axis. In this case LM-MNM algorithm may require more, than 1500 iterations to achieve performance comparable to CMNM. According to the averaged learning curve in Figures~\ref{fig:cvcnn_complex} and~\ref{fig:cvcnn_imag} LM-MNM and CMNM converge faster compared to LM-NM and CNM, which are based on a full Hessian computation.

\begin{figure}[ht]
  \centering
  \begin{minipage}[b]{0.48\textwidth}
    \centering
    \includegraphics[width=\textwidth]{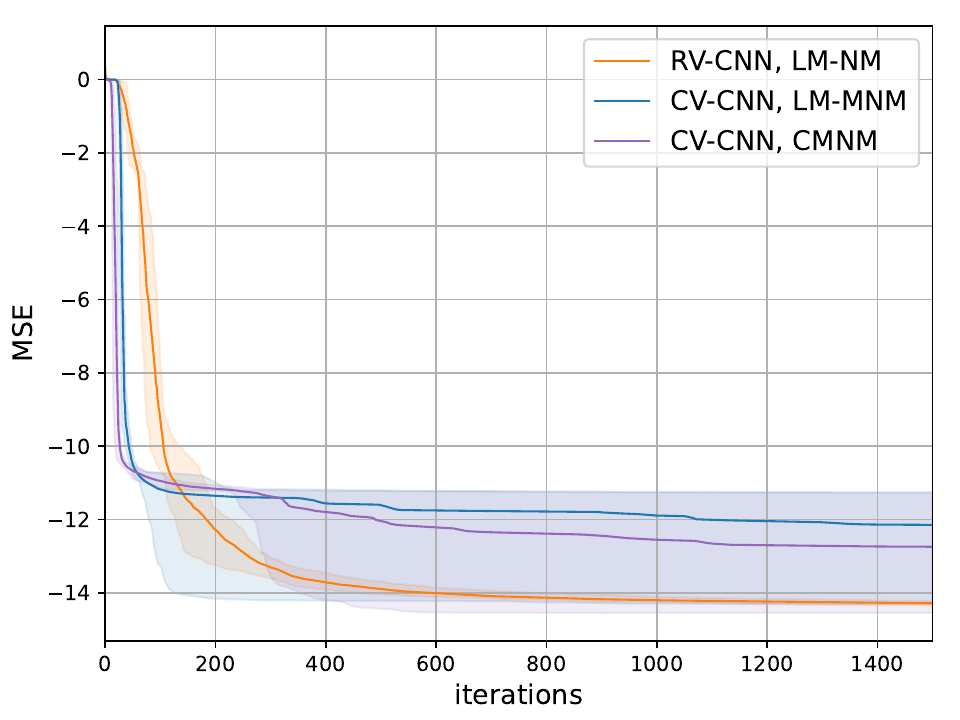}
  \end{minipage}
  \begin{minipage}[b]{0.48\textwidth}
    \centering
    \includegraphics[width=\textwidth]{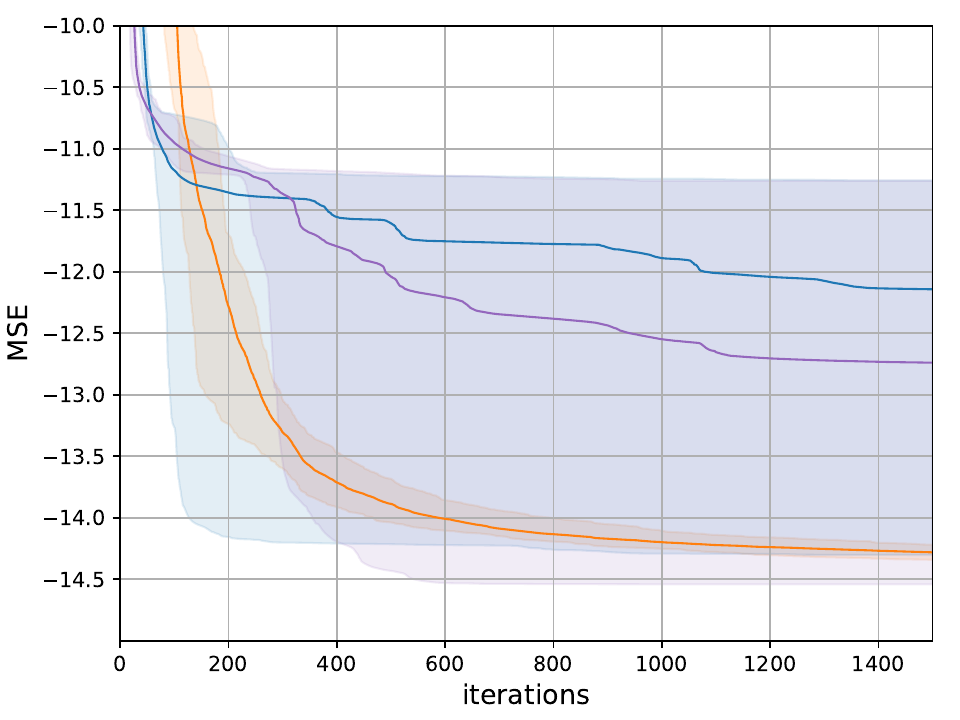}
  \end{minipage}
  \captionsetup{justification=centering}
  \caption{Comparison of learning curves for RV-CNN trained by LM-NM against CV-CNN trained by LM-MNM and CMNM.}
  \label{fig:compare}
\end{figure}

In addition, on Figure~\ref{fig:compare} we compare the convergence speed of RV-CNN, optimized by LM-NM, with that of CV-CNN, trained by LM-MNM and CMNM. One observes that RV-CNN in average converges in cost value to NMSE~=~-14.28~dB, which is better than the values NMSE~=~-12.76~dB and NMSE~=~-12.19~dB for CMNM and LM-MNM, correspondingly. Nevertheless, the best performance which was achieved by RV-CNN and CV-CNN when trained with the mentioned methods is approximately the same $\approx$ -14.5~dB (Table~\ref{tab:converge_eval}).

Furthermore, Table~\ref{tab:converge_eval} shows that LM-MNM~\eqref{MNMiterate} and CMNM~\eqref{RMNM}, which are based on use of the mixed Hessian~\eqref{MNMiterate} require $\approx$ 3.5 times less time per iteration comparing to LM-NM for RV-CNN and $\approx$~5 times less comparing to LM-NM for CV-CNN. This may be explained by the computation of the mixed Hessian using the Jacobian of the model output w.r.t. the model parameters, as shown in Theorem~\ref{thm:main}.

Finally, CV-CNN requires $\approx$ 1.5 times less real parameters compared to RV-CNN (Table~\ref{tab:table_telecom_hyper}) to achieve the same performance when trained by LM-NM (Table~\ref{tab:converge_eval}).

\begin{table}[ht]
\centering
\caption{Comparison of performance, relative convergence time of the proposed methods.}
\begin{tabular}{|c|c|c|c|c|c|c|}
\hline
    & \begin{tabular}[c]{@{}c@{}}RV-CNN,\\ LM-NM\end{tabular} & \begin{tabular}[c]{@{}c@{}}RV-CNN,\\ CNM\end{tabular} & \begin{tabular}[c]{@{}c@{}}CV-CNN,\\ LM-NM\end{tabular} & \begin{tabular}[c]{@{}c@{}}CV-CNN,\\ CNM\end{tabular} & \begin{tabular}[c]{@{}c@{}}CV-CNN,\\ LM-MNM\end{tabular} & \begin{tabular}[c]{@{}c@{}}CV-CNN,\\ CMNM\end{tabular} \\ \hline
    \begin{tabular}[c]{@{}c@{}}Relative \\ average time  \\ per iter.\end{tabular} & 1 & 0.98 & 1.39 & 1.44 & \textbf{0.28} & 0.30 \\ \hline
    \begin{tabular}[c]{@{}c@{}}Average\\ NMSE, \\ dB\end{tabular} & \textbf{-14.28} & -13.62 & -11.91     & -12.00 & -12.19 & -12.76 \\ \hline
    \begin{tabular}[c]{@{}c@{}}Best \\ NMSE,\\  dB\end{tabular} & -14.34 & -13.91 & -14.23     & -14.12 & -14.30 & \textbf{-14.55} \\ \hline
\end{tabular}
\label{tab:converge_eval}
\end{table}

\subsubsection{Simulations for the Abalone task}

\begin{figure}[ht]
  \centering
  \begin{subfigure}[t]{0.48\textwidth}
    \centering
    \includegraphics[width=\textwidth]{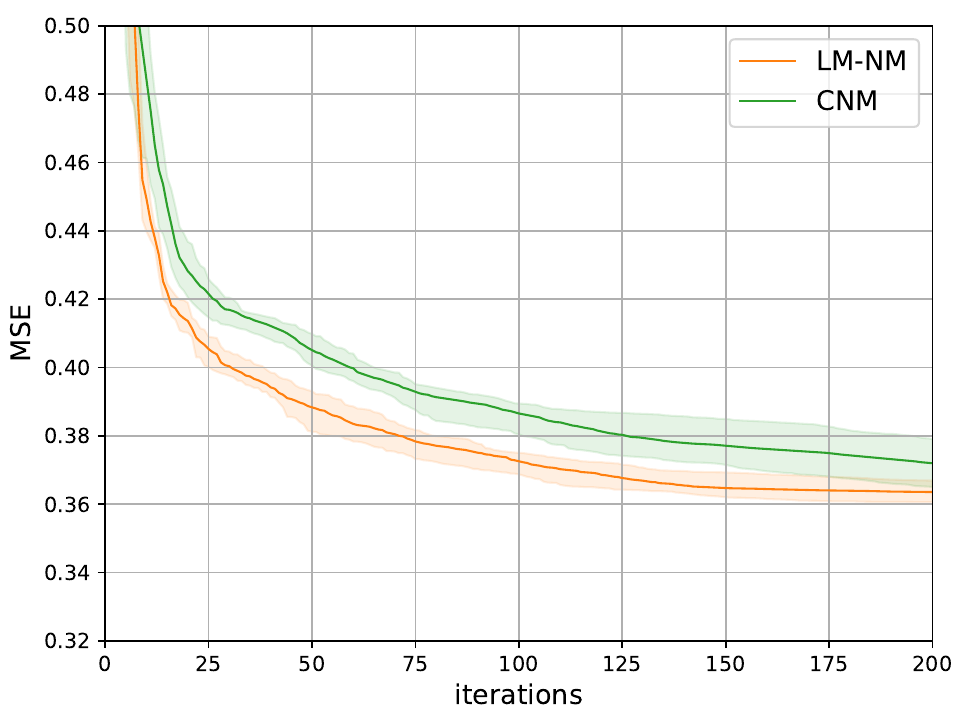}
    \captionsetup{justification=centering}
    \caption{Real weights.}
    \label{fig:abalone_real}
  \end{subfigure}
  \hfill
  \begin{subfigure}[t]{0.48\textwidth}
    \centering
    \includegraphics[width=\textwidth]{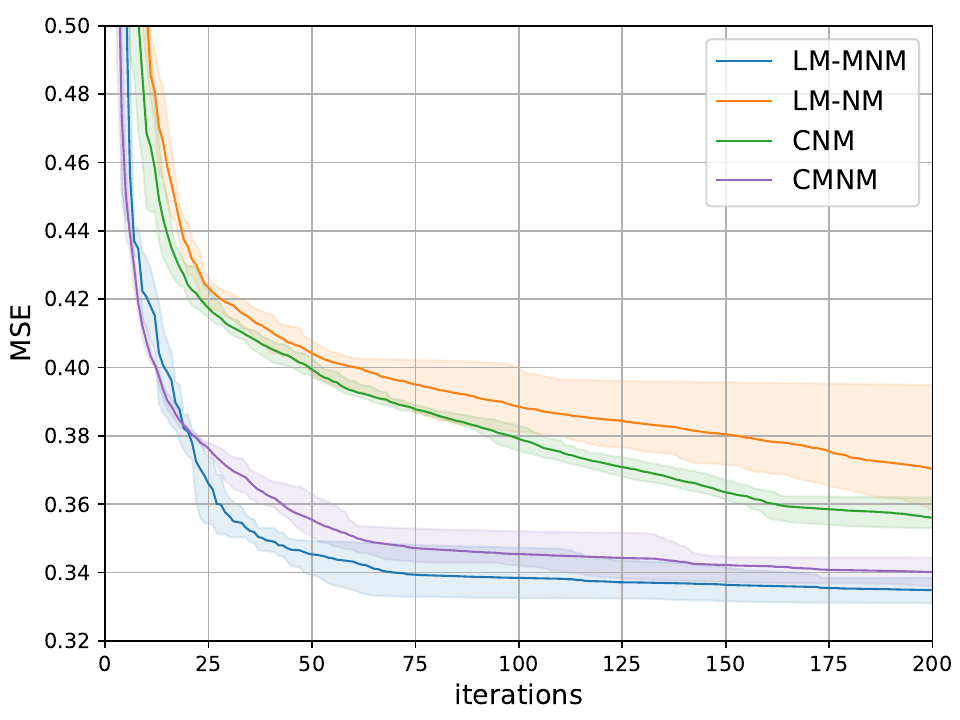}
    \captionsetup{justification=centering}
    \caption{Complex weights. Initial parameters in the complex plane.}
    \label{fig:abalone_complex}
  \end{subfigure}
  \\
  \begin{subfigure}[t]{0.48\textwidth}
    \centering
    \includegraphics[width=\textwidth]{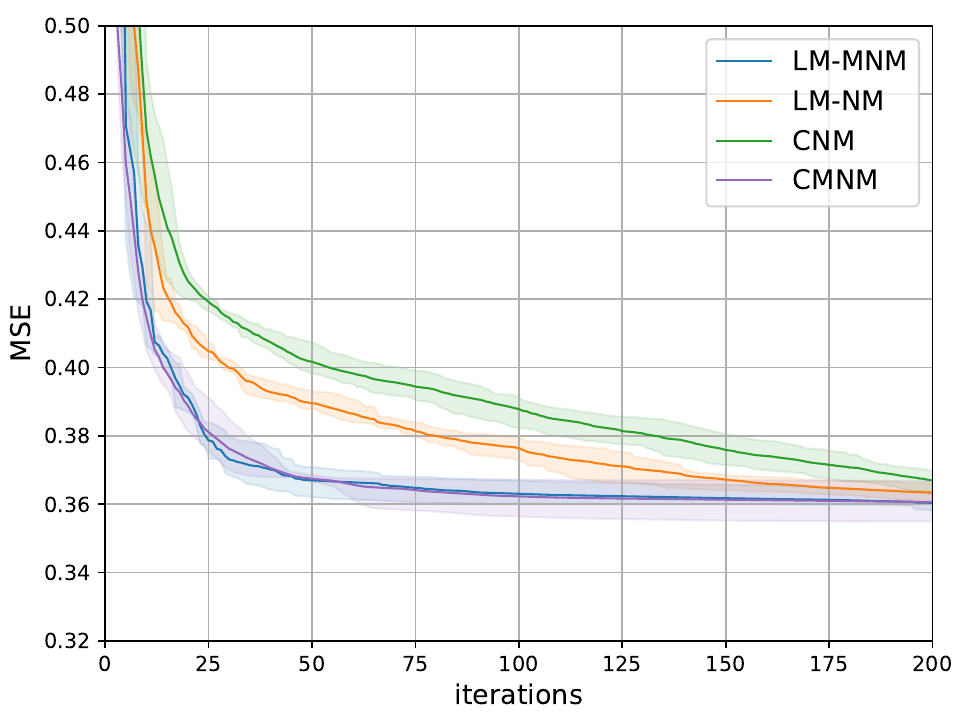}
    \captionsetup{justification=centering}
    \caption{Complex weights. Initial parameters on the real axis.}
    \label{fig:real_init}
  \end{subfigure}
  \hfill
  \begin{subfigure}[t]{0.48\textwidth}
    \centering
    \includegraphics[width=\textwidth]{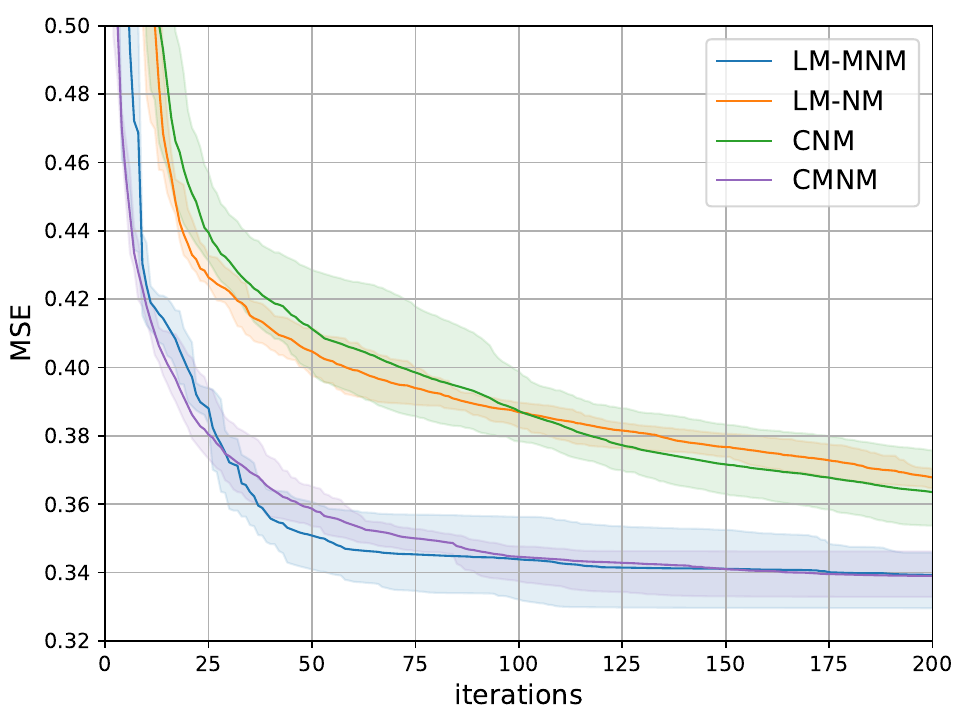}
    \captionsetup{justification=centering}
    \caption{Complex weights. Initial parameters on the imaginary axis.}
    \label{fig:imaginary_init}
  \end{subfigure}
  \captionsetup{justification=centering}
  \caption{Learning curves for the tested optimization algorithms using regression NN models with different weight configurations.}
  \label{fig:abalone}
\end{figure}

Abalone is a classical machine learning problem aimed at predicting the age of abalone from physical measurements. The dataset consists of 4177 objects with $M=8$ continuous normalized features and a continuous normalized target. Based on this data, we trained a small feed-forward fully connected neural network (NN) with one hidden layer of $H=10$ neurons, using tanh as activation function and mean squared error (MSE) as loss function. All features and targets in the dataset are real-valued, however we extended the problem to the complex space in three out of four cases to potentially improve convergence. The total number of parameters is $(M+1)H+(H+1)=101$. In each experiment the initial weights were selected from a normal distribution with mean $= 0$, standard deviation $= 0.1$, either along the real or imaginary axis, or in the entire complex space. The results are presented in Figure~\ref{fig:abalone}.

\begin{table}[ht]
\centering
\caption{Performance of the tested methods applied to the abalone problem.}
\begin{tabular}{|c|c|c|c|c|c|c|}
\hline
    & \begin{tabular}[c]{@{}c@{}}Real,\\ LM-NM\end{tabular} & \begin{tabular}[c]{@{}c@{}}Real,\\ CNM\end{tabular} & \begin{tabular}[c]{@{}c@{}}Complex,\\ LM-NM\end{tabular} & \begin{tabular}[c]{@{}c@{}}Complex,\\ CNM\end{tabular} & \begin{tabular}[c]{@{}c@{}}Complex,\\ LM-MNM\end{tabular} & \begin{tabular}[c]{@{}c@{}}Complex,\\ CMNM\end{tabular} \\ \hline
    \begin{tabular}[c]{@{}c@{}}Average\\ MSE \end{tabular} & 0.364 & 0.372 & 0.370     & 0.356 & \textbf{0.334} & 0.340 \\ \hline
    \begin{tabular}[c]{@{}c@{}}Best \\ MSE \end{tabular} & 0.361 & 0.365 & 0.358     & 0.353 & \textbf{0.331} & 0.336 \\ \hline
    \begin{tabular}[c]{@{}c@{}}Average\\ $R^2$ \end{tabular} & 0.636 & 0.628 & 0.630     & 0.644 & \textbf{0.665} & 0.660 \\ \hline
    \begin{tabular}[c]{@{}c@{}}Best \\ $R^2$ \end{tabular} & 0.639 & 0.635 & 0.642     & 0.647 & \textbf{0.669} & 0.664 \\ \hline
\end{tabular}
\label{tab:abalone}
\end{table}

\begin{figure}[ht]
  \centering
  \begin{minipage}[b]{0.48\textwidth}
    \centering
    \includegraphics[width=\textwidth]{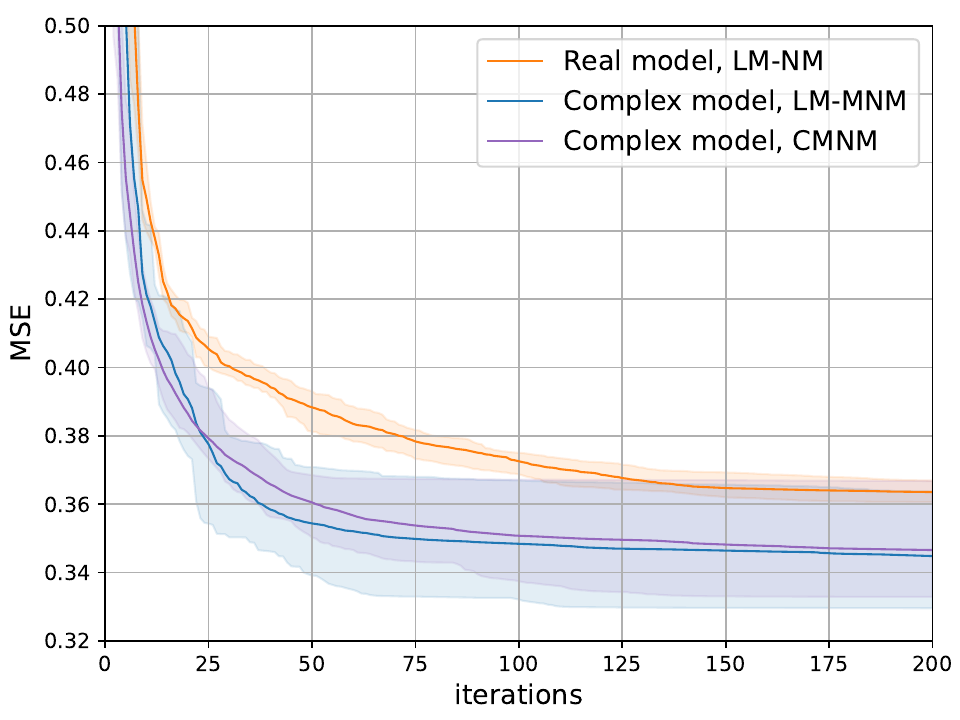}
  \end{minipage}
  \begin{minipage}[b]{0.48\textwidth}
    \centering
    \includegraphics[width=\textwidth]{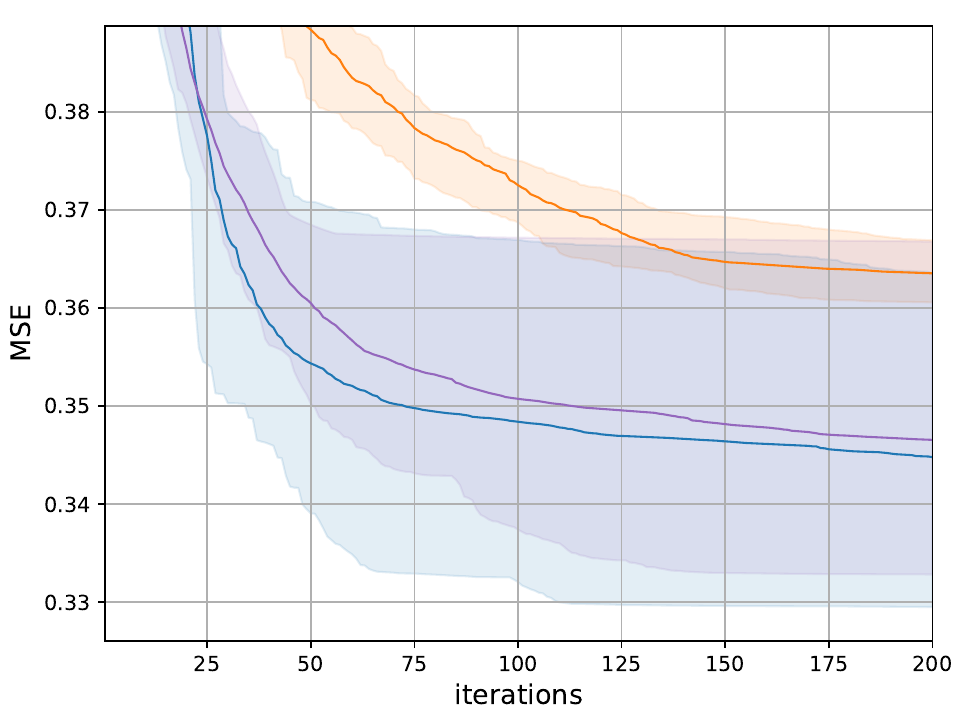}
  \end{minipage}
  \captionsetup{justification=centering}
  \caption{Comparison of learning curves for training of models with real and complex parameters.}
  \label{fig:compare_abalone}
\end{figure}

In terms of complex-weight model training, it can be seen that the LM-MNM and CMNM methods showed a higher convergence rate, approximately four times faster than both LM-NM and CNM: ~50 iterations versus ~200 iterations in all cases. Figure~\ref{fig:real_init} shows that when starting from points on the real axis, the model does not achieve the same performance as when starting from the entire complex plane or the imaginary axis, as shown in Figures~\ref{fig:abalone_complex} and~\ref{fig:imaginary_init} (MSE $\approx0.36$ vs $0.34)$. This is because, in this case, the weights do not leave the real line, effectively halving the number of trained parameters, as in the real model. Nevertheless, models with complex weights exploitation allows using efficient LM-MNM and CMNM methods, which provide $\approx4$ times convergence rate improvement and $0.02$ MSE improvement (after 200 iterations) in comparison with real-weight model, trained by LM-NM, which is shown in~Figure~\ref{fig:compare_abalone}.

The summary of experiments depicted in Figures~\ref{fig:abalone_real} and~\ref{fig:abalone_complex} including the $R^2$-score as a quality metric are presented in Table~\ref{tab:abalone}. All methods use a similar time per iteration, but the convergence of LM-MNM and CMNM is faster due to the smaller number of iterations. Moreover, LM-MNM achieves the best performance $R^2=0.699$ compared to LM-NM with $R^2=0.642$ for the same model.

\end{document}